\documentclass[a4paper, authoryear]{article}
\usepackage{techRepsPDF}
\usepackage[parfill]{parskip}
\usepackage{natbib}
\usepackage[ruled,vlined,linesnumbered]{myalgorithm2e}
\SetCommentSty{textsf}
\usepackage{graphicx} 
\usepackage{placeins}
\usepackage{newproof}
\usepackage{setspace}
\usepackage{amsmath, amssymb}
\usepackage{threeparttable, longtable}
\usepackage[margin=1cm]{caption}
\usepackage[pdftitle={Multiple-choice Vector Bin Packing: Arc-flow Formulation}, pdfauthor={Filipe Brandao}]{hyperref}
\usepackage{fullpage}
\newtheorem{example}{Example}

\def\HD#1#2{\vrule height #1pt depth #2pt width 0pt\relax}
\def\up{\HD{10}{0}}
\def\down{\HD{0}{5}}

\newcommand{\sip}{\mathop{\mathrm{ip}}}
\newcommand{\tip}{t^{\sip}}

\newcommand{\sbb}{\mathop{\mathrm{bb}}} 
\newcommand{\nbb}{n^{\sbb}}

\newcommand{\stot}{\mathop{\mathrm{tot}}} 
\newcommand{\ttot}{t^{\stot}}

\newcommand{\vS}{\textsc{s}}
\newcommand{\vT}{\textsc{t}}

\def\ptitleA{
Multiple-choice Vector Bin Packing:\\
Arc-flow Formulation with\\ Graph Compression
}
\def\ptitleB{
Multiple-choice Vector Bin Packing:\\
Arc-flow Formulation with Graph Compression
}
\def\pnumber{DCC-2013-13}
\def\pauthor{
Filipe Brand\~ao\\
{\small INESC TEC and Faculdade de Ci\^{e}ncias,
Universidade do Porto, Portugal}\\
{\small\texttt{fdabrandao@dcc.fc.up.pt}}
\ \\ 
\ \\
Jo\~ao Pedro Pedroso\\
{\small INESC TEC and Faculdade de Ci\^{e}ncias,
Universidade do Porto, Portugal}\\
{\small\texttt{jpp@fc.up.pt}}
}

\graphicspath{{figures_flow/}}

\begin{document}

\trtitle{\ptitleA}
\trauthor{\pauthor}
\trnumber{\pnumber}
\mkcoverpage

\title{\textbf{\ptitleB}}
\author{\pauthor}
%\date{\today}
\date{December 12, 2013}
\maketitle

\begin{abstract}
The vector bin packing problem (VBP) is a generalization of bin packing with multiple constraints. 
In this problem we are required to pack items, represented by $p$-dimensional vectors, 
into as few bins as possible. 
The multiple-choice vector bin packing (MVBP)
is a variant of the VBP
in which bins have several types and items have several incarnations.
We present an exact method, based on an arc-flow formulation
with graph compression, for solving MVBP
by simply representing all the patterns in a very compact graph.
As a proof of concept we report computational results
on a variable-sized bin packing data set.
\end{abstract}

%\doublespacing
%\singlespacing
%\onehalfspacing

\section{Introduction}

The vector bin packing problem (VBP), 
also called general assignment problem by
some authors, is a generalization of bin packing with multiple constraints. 
In this problem, we are
required to pack $n$ items of $m$ different types, represented by $p$-dimensional vectors, 
into as few bins as possible.
The multiple-choice vector bin packing problem (MVBP) is a variant 
of VBP in which bins have several types (i.e., sizes and costs) 
and items have several incarnations (i.e., will take one of several possible sizes);
this occurs typically in situations where one of several incompatible 
decisions has to be taken (see, e.g., \citealt{Patt-Shamir:2012:VBP:2206436.2206615}).

\cite{BrandaoGeneralArcFlow} present a general arc-flow formulation
with graph compression for vector packing.
This formulation is equivalent to the model of \cite{gomory2},
thus providing a very strong linear relaxation. 
It has proven to be very effective
on a large variety of problems through reductions to vector packing.
In this paper, we apply the general arc-flow formulation 
to the multiple-choice vector packing problem.

The remainder of this paper is organized as follows.
Section~\ref{sec:arcflow} presents the arc-flow formulation for MVBP.
Some computational results are presented in Section~\ref{sec:results}
and Section~\ref{sec:conclusions} presents the conclusions.

\section{Arc-flow formulation with graph compression for MVBP}
\label{sec:arcflow}

In order to solve a cutting/packing problem,
the arc-flow formulation proposed in \cite{BrandaoGeneralArcFlow}
only requires the corresponding directed acyclic multigraph $G=(V,A)$ 
containing every valid packing pattern represented as a path from
the source to the target.
In order to model MVBP, we will start by defining the underlying graph.

For a given $i$, let $\mathbf{J}_i$ be the set of incarnations of item $i$, and
let $\mathbf{I} = \{(i,j) : i=1..m, j \in \mathbf{J}_i\}$
be the set of items.
Let $it_i^j=(i,j)\in \mathbf{I}$ be the incarnation~$j$ of item~$i$
and $w(it_i^j)$ its weight vector.
For the sake of simplicity, we define~$it_0^0$
as an item with weight zero in every dimension; this artificial item
is used to label loss arcs.
Let~$b_i$ be the demand of items of type~$i$, for $i=1,\ldots,m$.
Let~$q$ be the number of bin types.
Let~$W(t)$ and~$C(t)$ be the capacity vector
and the cost of bins of type~$t$, respectively.

\begin{figure}[h!tbp]
\caption{Arc-flow graph containing every valid packing pattern
for Example~\ref{ex:example1}.\label{fig:graph1}}

  \centering
  \includegraphics[scale=1]{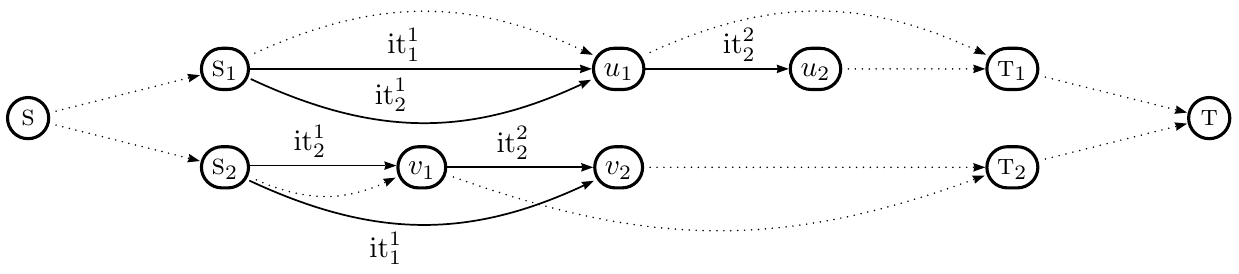}

\caption*{\footnotesize 
Paths from $\vS_t$ to $\vT_t$
represent every valid pattern for bins of type $t$,
for each $t$.
Paths from $\vS$ to $\vT$ represent 
every valid packing pattern using any bin type.
}
\end{figure}

\begin{example}
\label{ex:example1}
Figure~\ref{fig:graph1} shows the graph associated with a
two dimensional ($p=2$) instance with bins of two types ($q=2$).
The bins of type 1 have capacity $W(1) = (100,75)$ and cost $C(1)=3$.
The bins of type 2 have capacity $W(2) = (75,50)$ and cost $C(2)=2$.
There are three items ($n=3$) to pack of two different types ($m=2$).
The first item type has demand $b_1=2$, and a single incarnation with weight $w(it_1^1)=(75, 50)$.
The second item type has demand $b_2=1$, and two incarnations with weights $w(it_2^1) =(40, 15)$
and $w(it_2^2)=(25, 25)$.
\end{example}

We need to build a graph for each bin type considering every item incarnation 
as a different item. 
The arc-flow graphs must contain every valid packing
pattern represented as a path from the source to the target,
and they may not contain any invalid pattern.
These graphs -- say, $G_1,\ldots,G_t$ -- can be built using the step-by-step
algorithm proposed in~\cite{MThesisBrandao} or the
algorithm proposed in~\cite{BrandaoGeneralArcFlow} 
(recommended for efficiency).
Both algorithms perform graph compression
and hence the resulting graphs tend to be small.

Figure~\ref{fig:graph1} shows an arc-flow graph 
for Example~\ref{ex:example1}.
Paths from $\vS_t$ to $\vT_t$
represent every valid pattern for bins of type $t$,
for $t=1,\ldots,q$. 
Each of these subgraphs is built considering
every item incarnation as a different item.
We connect a super source node $\vS$ to every $\vS_t$,
and every $\vT_t$ to a super target node $\vT$.
Paths from $\vS$ to $\vT$ represent 
every valid packing pattern using any bin type.

Graphs $G_1,\ldots,G_t$ are already compressed, but
in order to reduce the whole graph size even more, 
we apply again to $G$ the final compression step of
the method proposed in~\cite{BrandaoGeneralArcFlow}.
Note that this compression step can only be applied
if the set of item incarnations does not depend on the bin type.
We relabel the graph using the longest paths from the source in each dimension.
Let ($\psi^1(v)$, $\psi^2(v)$, \ldots, $\psi^p(v)$) be the label of node~$v$ in the
final graph, where
\begin{alignat}{3}
\psi^d(v) & = & \left\{ \begin{array}{ll}
                0 & \mbox{if }  v = \vS, \\
                \max_{(u,v',it_i^j) \in A:v'=v}\{\psi^{d}(u) + w(it_i^j)^d\}  & \mbox{otherwise.}\\
                \end{array}\right.           
\end{alignat}

\begin{figure}[h!tbp]
\caption{After applying the final compression step to the graph of Figure~\ref{fig:graph1}.\label{fig:graph2}}

  \centering
  \includegraphics[scale=1]{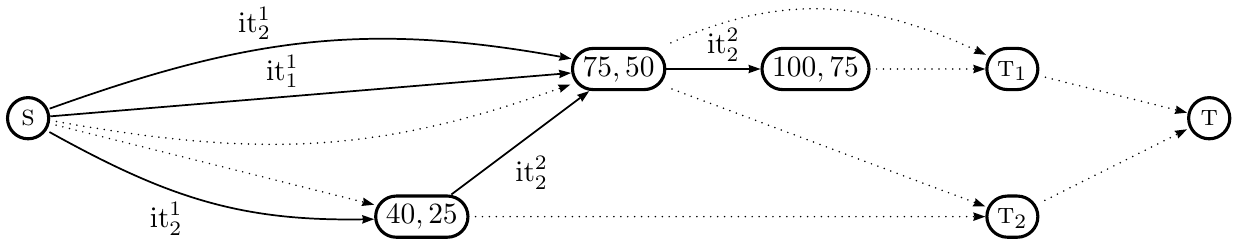}

\caption*{\footnotesize 
The final compression step removed 3 vertices and 3 arcs on this small example.}
\end{figure}

The final graph may contain parallel arcs for different incarnations of the same item.
Since having multiple parallel arcs for the same item is redundant, only one of them
is left.

The arc-flow formulation for multiple-choice vector bin packing is the following:
\begin{alignat}{3}
  & \mbox{minimize }   && \sum_{t=0}^q C(t) f_{\vT_t,\vT,it_0^0}\label{eq:new1}\\  
  & \mbox{subject to } \qquad && \sum_{(u,v,it)\in A:v=k} \hspace{-5mm}f_{u,v,it} \hspace{3mm}-\hspace{-3mm} \sum_{(v,r,it) \in A:v=k} \hspace{-5mm}f_{v,r,it} = &&
  \left\{ 
    \begin{array}{rl}
      -z & \mbox{if } k = \vS, \\
      z & \mbox{if } k = \vT,\\
      0 & \mbox{for } k \in V \setminus \{\vS, \vT \},\\
    \end{array} \label{eq:new2}
  \right.\\    
      &&         & \sum_{(u, v, it_i^j) \in A:i=k} \hspace{-5mm}f_{u,v,it_i^j} \geq b_k, && k \in \{1,\ldots,m\} \setminus J, \label{eq:new3}\\ 
      &&         & \sum_{(u, v, it_i^j) \in A:i=k} \hspace{-5mm}f_{u,v,it_i^j} = b_k, && k \in J,  \label{eq:new4}\\
      &&         & f_{u,v,it_i^j} \leq b_i, && \forall (u,v,it_i^j) \in A, \mbox{ if } i \neq 0, \label{eq:new5}\\
      &&         & f_{u,v,it_i^j} \geq 0, \mbox{ integer}, && \forall (u,v,it_i^j) \in A, \label{eq:new6}      
\end{alignat}
where $z$~can be seen as a feedback from $\vT$ to $\vS$;
$m$~is the number of different items;
$q$~is the number of bin types;
$b_i$~is the demand of items of type $i$;
$V$~is the set of vertices,
$\vS$~is the source vertex and $\vT$~is the target;
$A$~is the set of arcs, where
each arc has three components~$(u, v, it_i^j)$
corresponding to an arc between nodes~$u$ and~$v$
that contributes to the demand of items of type $i$;
arcs~$(u,v,it_0^0)$ are loss arcs;
$f_{u,v,it}$~is the amount of flow along the arc~$(u, v, it)$;
and $J \subseteq \{1,\ldots,m\}$ is a subset of items whose demands 
are required to be satisfied exactly for efficiency purposes.
For having tighter constraints, one may set $J = \{i = 1,\ldots, m : b_i = 1\}$ 
(we have done this in our experiments). The main difference
between this and the original arc-flow formulation is the objective function.

Algorithm~\ref{alg:solve} illustrates our solution method.
More details on algorithms for graph construction and solution extraction
are given in~\cite{BrandaoGeneralArcFlow} and~\cite{MThesisBrandao}.

\begin{algorithm}[!h]
\caption{MVBP Solution Method}\label{alg:solve}
\SetKwInOut{Input}{input}
\SetKwInOut{Output}{output}

\SetKwFunction{buildGraph}{buildGraph}
\SetKwFunction{solve}{solveMVBP}
\SetKwFunction{compress}{compress}
\SetKwFunction{MIPSolver}{MIPSolver}
\SetKwFunction{extractSol}{extractSolution}
\SetKwData{weight}{weight}
\SetKwData{demand}{demand}
\SetKwData{xlabels}{labels}
\SetKwData{varmx}{mx}
\SetKwData{varmi}{mi}
\SetKwData{NIL}{NIL}

\SetKwBlock{Function}{}{}

\Input{
$\mathbf{I}$~-~set of items; 
$m$~-~number of different items;
$w(it_i^j)$~-~weight vector of the incarnation~$j$ of item~$i$; 
$b_i$~-~demand of item~$i$; 
$q$~-~number of bin types;
$W(t)$, $C(t)$~-~capacity and cost of bins of type~$t$, respectively; }
\Output{MVBP Solution}

\textbf{function} $\solve(\mathbf{I}, m, w, b, q, W, C)$:
\Function{
$V, A \gets (\{\vS, \vT\},\emptyset)$\;
\For(\tcp*[f]{for each bin type $t$}){$t\leftarrow 1$ \KwTo $q$}{       
    $\xlabels \gets \mathbf{I}$\;
    $\weight \gets [w(it_i^j) : it_i^j \in \mathbf{I}]$\;
    $\demand \gets [b_i : it_i^j \in \mathbf{I}]$\;
    $(G_t, \vS_t, \vT_t) \gets \buildGraph(m, \xlabels, \weight, \demand, W(t))$\tcp*{build the arc-flow graph $G_t$ for bins of type $t$}
    $(V_t, A_t) \gets G_t$\;
    $V \gets V \cup V_t$\;
    $A \gets A \cup A_t \cup \{(\vS, \vS_t, it_0^0), (\vT_t, \vT, it_0^0)\}$\;
}
$G \gets (V, A)$\;
$G \gets \compress(G)$\tcp*{apply the final compression step to $G$}
$f \gets \MIPSolver(\mbox{arc-flow model}, q, C, G, b)$\tcp*{solve the arc-flow model over $G$}
\Return $\extractSol(f)$\tcp*{extract the MVBP solution from the arc-flow solution}
}
\end{algorithm}

\section{Computational results}
\label{sec:results}

As a proof of concept, we used to the arc-flow formulation
to solve variable-sized bin packing.
A specific method for solving this problem has been presented in~\cite{Alves20071333};
here we proposed to simply model it as a unidimensional multiple-choice
vector bin packing problem.
We solved the benchmark data set of~\cite{MonaciThesis}, which is composed of 300 instances.
In this data set, the items sizes were randomly generated 
within three different ranges: $w_i \in [1,100]$ (X=1), $w_i \in [20,100]$ (X=2) 
and $w_i \in [50,100]$ (X=3).
There are instances with three ($q=3$) and five $(q=5)$ bin types;
the bin sizes are $[100,120,150]$ for $q=3$ and $[60,80,100,120,150]$ for $q=5$.
For each range and each number of bin types,
there are 10 instances for each $n \in \{25, 50, 100, 200, 500\}$.
The average run time on the 300 instances was less than 1 second
and none of these instances took longer than 6 seconds to be solved exactly.

\pagebreak

Table~\ref{tab:results} presents the results.
The meaning of each column is as follows: 
$q$~-~number of different bin types;
$n$~-~ total number of items;
$m$~-~number of different items;
$\#v$, $\#a$~-~number of vertices and arcs in the final arc-flow graph;
$\%v, \%a$~-~percentage of vertices and arcs removed by the final compression step;
$\tip$~-~time spent solving the model; 
$\nbb$~-~average number of nodes explored in the branch-and-bound procedure; 
$\ttot$~-~run time in seconds.
The values shown are averages over the 10 instances in each class. 

CPU times were obtained using a computer with two Quad-Core Intel 
Xeon at 2.66GHz, running Mac OS X 10.8.5.
The graphs for each bin type were generated using the algorithm
proposed in~\cite{BrandaoGeneralArcFlow}, which was implemented in \texttt{C++},
and the final model was produced using \texttt{Python}.
The models were solved using 
\texttt{Gurobi} 5.5.0 (\citealt{gurobi}), a state-of-the-art mixed integer programming solver.
The parameters used in \texttt{Gurobi} were Threads~=~1 (single
thread), Presolve~=~1 (conservative), Method~=~2 (interior point methods),
MIPFocus~=~1 (feasible solutions), Heuristics~=~1, MIPGap~=~0, MIPGapAbs~=~$1-10^{-5}$ 
and the remaining parameters were \texttt{Gurobi}'s default values. 
The branch-and-cut solver used in \texttt{Gurobi} uses a series of cuts; in our
models, the most frequently used were Gomory, Zero half and MIR. %Zero half: 72, MIR: 41, Gomory: 90
The source code is available online\footnote{\url{http://www.dcc.fc.up.pt/\~fdabrandao/code}}.

\begin{table}[h!tbp]
\centering\begin{threeparttable}[b]
\scriptsize 
\caption{Results for variable-sized bin packing.\label{tab:results}}
\begin{tabular}{rrrrrrrrrrr}
\hline
\up\down
Range & $q$ & $n$ & $m$ & $\#v$ & $\#a$ & $\%v$ & $\%a$ & $\nbb$ & $\tip$ & $\ttot$\\
\hline	
\up
X=1 & 3 & 25 & 22.0 & 71.1 & 616.9 & 41.90 & 13.68 & 13.5 & 0.15 & 0.23\\
X=1 & 3 & 50 & 38.3 & 115.1 & 1,641.7 & 53.38 & 26.39 & 7.2 & 0.56 & 0.69\\
X=1 & 3 & 100 & 62.3 & 134.3 & 3,128.8 & 56.89 & 38.49 & 3.2 & 1.70 & 1.94\\
X=1 & 3 & 200 & 86.8 & 142.9 & 4,930.3 & 57.86 & 43.33 & 0.0 & 1.35 & 1.72\\
X=1 & 3 & 500 & 98.8 & 148.3 & 6,047.9 & 58.34 & 46.07 & 0.0 & 1.41 & 1.91\\
\up
X=1 & 5 & 25 & 22.4 & 91.8 & 807.8 & 48.64 & 16.40 & 0.0 & 0.11 & 0.21\\
X=1 & 5 & 50 & 39.1 & 122.9 & 1,947.6 & 60.62 & 27.83 & 0.8 & 0.57 & 0.73\\
X=1 & 5 & 100 & 61.8 & 139.3 & 3,518.9 & 66.12 & 40.94 & 0.0 & 0.73 & 1.01\\
X=1 & 5 & 200 & 85.6 & 147.1 & 5,289.6 & 68.06 & 48.89 & 0.0 & 2.19 & 2.62\\
X=1 & 5 & 500 & 98.5 & 151.9 & 6,459.6 & 69.02 & 52.91 & 0.0 & 1.26 & 1.85\\
\up
X=2 & 3 & 25 & 21.6 & 38.4 & 307.0 & 35.46 & 10.27 & 0.0 & 0.04 & 0.10\\
X=2 & 3 & 50 & 36.9 & 72.8 & 1,020.7 & 43.19 & 10.25 & 0.0 & 0.12 & 0.21\\
X=2 & 3 & 100 & 58.7 & 96.0 & 2,379.0 & 50.36 & 13.30 & 0.0 & 0.33 & 0.49\\
X=2 & 3 & 200 & 73.1 & 108.7 & 3,554.2 & 53.70 & 16.11 & 0.0 & 0.67 & 0.90\\
X=2 & 3 & 500 & 80.0 & 113.3 & 4,076.2 & 54.39 & 17.39 & 0.0 & 0.41 & 0.68\\
\up
X=2 & 5 & 25 & 22.2 & 43.3 & 348.7 & 39.33 & 11.63 & 0.0 & 0.02 & 0.09\\
X=2 & 5 & 50 & 37.4 & 75.3 & 1,037.8 & 46.88 & 10.82 & 0.0 & 0.09 & 0.20\\
X=2 & 5 & 100 & 57.5 & 96.0 & 2,275.6 & 54.26 & 13.10 & 0.0 & 0.30 & 0.47\\
X=2 & 5 & 200 & 73.7 & 110.2 & 3,757.8 & 59.99 & 17.17 & 0.0 & 0.49 & 0.75\\
X=2 & 5 & 500 & 79.7 & 115.5 & 4,267.1 & 60.94 & 18.75 & 0.0 & 0.84 & 1.20\\
\up
X=3 & 3 & 25 & 19.4 & 15.6 & 90.4 & 16.69 & 3.48 & 0.0 & 0.00 & 0.06\\
X=3 & 3 & 50 & 30.7 & 22.8 & 148.0 & 12.69 & 2.78 & 0.0 & 0.01 & 0.07\\
X=3 & 3 & 100 & 44.8 & 36.5 & 228.0 & 7.86 & 1.35 & 0.0 & 0.01 & 0.07\\
X=3 & 3 & 200 & 49.3 & 41.9 & 254.8 & 6.69 & 1.16 & 0.8 & 0.01 & 0.08\\
X=3 & 3 & 500 & 50.0 & 43.0 & 260.0 & 6.52 & 1.14 & 0.0 & 0.01 & 0.07\\
\up
X=3 & 5 & 25 & 18.5 & 18.2 & 101.0 & 25.03 & 6.82 & 0.0 & 0.00 & 0.06\\
X=3 & 5 & 50 & 31.6 & 26.3 & 179.8 & 20.14 & 4.83 & 0.0 & 0.00 & 0.08\\
X=3 & 5 & 100 & 43.3 & 36.5 & 251.1 & 15.62 & 3.23 & 0.0 & 0.01 & 0.08\\
X=3 & 5 & 200 & 49.5 & 44.4 & 297.9 & 13.62 & 2.65 & 0.0 & 0.01 & 0.09\\
\down
X=3 & 5 & 500 & 50.0 & 45.0 & 301.0 & 13.46 & 2.59 & 0.0 & 0.01 & 0.12\\
\hline
\end{tabular}
\end{threeparttable}
\end{table}

\section{Conclusions}
\label{sec:conclusions}

We propose an arc-flow formulation with graph compression
for solving multiple-choice vector bin packing problems.
This formulation is simple and proved to be effective
for solving variable-sized bin packing
as a unidimensional multiple-choice vector packing problem.
This paper shows the flexibility and effectiveness 
of the general arc-flow formulation with graph compression 
for modeling and solving cutting and packing problems,
beyond those solved through reductions to vector packing
in the original paper.

\FloatBarrier

% BibTeX users please use one of
%\bibliographystyle{plain}
%\bibliographystyle{spbasic}      % basic style, author-year citations
%\bibliographystyle{spmpsci}      % mathematics and physical sciences
%\bibliographystyle{spphys}       % APS-like style for physics
\bibliographystyle{apalike} 
\bibliography{paper}

\begin{thebibliography}{}

\bibitem[Alves and {Val\'erio de Carvalho}, 2007]{Alves20071333}
Alves, C. and {Val\'erio de Carvalho}, J. (2007).
\newblock Accelerating column generation for variable sized bin-packing
  problems.
\newblock {\em European Journal of Operational Research}, 183(3):1333 -- 1352.

\bibitem[Brand\~ao, 2012]{MThesisBrandao}
Brand\~ao, F. (2012).
\newblock {Bin Packing and Related Problems: Pattern-Based Approaches}.
\newblock Master's thesis, Faculdade de Ci\^encias da Universidade do Porto,
  Portugal.

\bibitem[Brand\~ao, 2013]{VPSolver}
Brand\~ao, F. (2013).
\newblock {VPSolver: Arc-flow Vector Packing Solver, Version 1.1}.
\newblock (Software program).

\bibitem[Brand\~ao and Pedroso, 2013]{BrandaoGeneralArcFlow}
Brand\~ao, F. and Pedroso, J.~P. (2013).
\newblock {Bin Packing and Related Problems: General Arc-flow Formulation with
  Graph Compression}.
\newblock Technical Report DCC-2013-08, Faculdade de Ci\^encias da Universidade
  do Porto, Portugal.

\bibitem[Gilmore and Gomory, 1963]{gomory2}
Gilmore, P. and Gomory, R. (1963).
\newblock {A linear programming approach to the cutting stock problem--part
  {II}}.
\newblock {\em Operations Research}, 11:863--888.

\bibitem[Gu et~al., 2013]{gurobi}
Gu, Z., Rothberg, E., and Bixby, R. (2013).
\newblock {Gurobi Optimizer, Version 5.5.0}.
\newblock (Software program).

\bibitem[Monaci, 2001]{MonaciThesis}
Monaci, M. (2001).
\newblock {\em Algorithms for Packing and Scheduling Problems}.
\newblock PhD thesis, Universit\`a di Bologna.

\bibitem[Patt-Shamir and Rawitz, 2012]{Patt-Shamir:2012:VBP:2206436.2206615}
Patt-Shamir, B. and Rawitz, D. (2012).
\newblock Vector bin packing with multiple-choice.
\newblock {\em Discrete Appl. Math.}, 160(10-11):1591--1600.

\end{thebibliography}
 
\end{document}